\documentclass{amsart}
\usepackage{amsfonts}
\usepackage{amsmath}
\usepackage{amssymb}

\newtheorem{theorem}{Theorem}

\newtheorem{corollary}[theorem]{Corollary}

\begin{document}
\title{Positively curved shrinking Ricci solitons are compact}
\author{Ovidiu Munteanu}
\email{ovidiu.munteanu@uconn.edu}
\address{Department of Mathematics, University of Connecticut, Storrs, CT
06268, USA}
\author{Jiaping Wang}
\email{jiaping@math.umn.edu}
\address{School of Mathematics, University of Minnesota, Minneapolis, MN
55455, USA}
\maketitle

\begin{abstract}
We show that a shrinking Ricci soliton with positive sectional curvature
must be compact. This extends a result of Perelman in dimension three and 
improves a result of Naber in dimension four, respectively.
\end{abstract}

\section{Introduction}
In studying the singularities of Ricci flows, Hamilton \cite{H} has introduced
the concept of Ricci solitons and signified their importance. The importance
was further demonstrated in the work of Perelman \cite{P1,P2}, where his
classification result of three dimensional shrinking gradient Ricci solitons
played a crucial role in the affirmative resolution of the Poincar\'{e}
conjecture. Recall that a complete manifold $(M,g)$ is a shrinking gradient Ricci
soliton if the equation

\begin{equation*}
\mathrm{Ric}+\mathrm{Hess}\left( f\right) =\frac{1}{2}\,g
\end{equation*}%
holds for some function $f.$ Here, $\mathrm{Ric}$ is the Ricci curvature of 
$\left( M,g\right) $ and $\mathrm{Hess}\left( f\right) $ the Hessian of $f.$

The importance of shrinking Ricci solitons can be partially seen through a conjecture attributed to Hamilton. The conjecture asserts that the blow-ups 
around a type-I singularity point of a Ricci flow (see \cite{CLN} for definition) always converge to (nontrivial) gradient shrinking Ricci solitons. Important progress toward
the conjecture was made in \cite{N, S}, where it was shown that blow-up limits must be
gradient shrinking Ricci solitons. The nontriviality issue, which was explicitly raised
by Cao \cite{C1}, was later taken up by Enders, M\"uller and Topping \cite{EMT}.

Two dimensional shrinking Ricci solitons have been classified by Hamilton \cite{H}. For the three dimensional case, Perelman \cite{P2} made the breakthrough and concluded that a noncollapsing shrinking gradient Ricci soliton with bounded curvature must be a quotient of the sphere $\mathbb{S}^{3},$ or $\mathbb{R}^{3},$ or $\mathbb{S}^{2}\times \mathbb{R}.$ The result was later shown to be true without any extra assumptions through the effort of \cite{CCZ, NW} by a different approach.
See also \cite{N, PW} for related works.

For a three dimensional shrinking gradient Ricci soliton with bounded curvature, according to Ivey \cite{I} and Hamilton \cite{H}, its sectional curvature is necessarily nonnegative (see \cite{C} for more general result). On the other hand, the strong maximum principle in \cite{H2} implies that the universal cover of such a soliton must split off a line if its curvature is not strictly positive. Together with the fact that a compact three dimensional gradient shrinking Ricci soliton is necessarily covered by the standard sphere \cite{H1, I} (see \cite{ELM} for a different proof), one sees that the classification result of Perelman can be deduced from the following theorem.

\begin{theorem} (Perelman \cite{P2})
Let $\left( M,g\right) $ be a three dimensional noncollapsing gradient shrinking Ricci soliton with positive and bounded sectional curvature. Then $\left(
M,g\right) $ must be compact.
\end{theorem}

A natural question, as has been raised by Cao \cite{C2}, is whether the preceding theorem actually holds true in full generality for all dimension. The main purpose of this short note is to provide a complete answer to this question. 

\begin{theorem}
Let $\left( M,g\right) $ be a gradient shrinking Ricci soliton with nonnegative sectional curvature and positive Ricci curvature. Then $\left(M,g\right) $ must be compact.
\end{theorem}

Notice that we do not require any extra assumptions such as the manifold is noncollapsing or the curvature is bounded. Recall that
a simply connected shrinking gradient Ricci soliton $M$ with 
nonnegative sectional curvature is necessarily of the form
$M=\mathbb{R}^k\times N$ with $N$ having positive Ricci curvature.
This was stated as Corollary 4 in \cite{PW} under an integrability 
assumption on the Ricci curvature. However,  according 
to \cite{MS}, this assumption automatically holds true for all shrinking gradient Ricci solitons. Combining this with Theorem 2, we arrived at the following conclusion.

\begin{corollary}
Let $\left( M,g,f\right) $ be an $n$ dimensional shrinking gradient Ricci soliton 
with nonnegative sectional curvature. Then $\left( M,g\right) $ must be
compact, or a quotient of $\mathbb{R}^{n}$ or of the product 
$\mathbb{R}^{k}$ $\times N^{n-k}$  with $1\leq k\leq n-2,$ where $N$ is 
a compact simply connected shrinking Ricci soliton of dimension $n-k$ with
positive Ricci curvature. 
\end{corollary}

Since a compact simply connected shrinking gradient Ricci soliton with nonnegative curvature operator must be the round sphere by B\"ohm and Wilking \cite{BW}, the 
following corollary is immediate. This also answers a question stated 
in \cite{CLN} on page 389.

\begin{corollary}
Let $\left( M,g,f\right) $ be an $n$ dimensional shrinking gradient Ricci soliton with nonnegative curvature operator. Then $\left( M,g\right) $ must be
a quotient of the sphere $\mathbb{S}^n,$ or  $\mathbb{R}^{n},$ or the product 
$\mathbb{R}^{k}$ $\times \mathbb{S}^{n-k}$  with $1\leq k\leq n-2.$
\end{corollary}

In \cite{N}, Naber was able to obtain a partial extension of Perelman's result
to four dimensional shrinking gradient Ricci solitons with positive curvature operator. 
As a consequence, he also proved Corollary 4  for the case $n=4$ under the further condition that the curvature is bounded.
However, in both \cite{P2} and \cite{N}, the proofs seem to be dimension specific
and it is not at all clear to us if they may be generalized to deal with arbitrary dimension case. Our proof is very much inspired by \cite{CLY}, where the authors obtained a 
quadratic decay positive lower bound for the scalar curvature of nontrivial, noncompact shrinking Ricci solitons. Here, we managed to obtain a similar estimate for the Ricci curvature, which in turn implies that the scalar curvature must increase and become unbounded at infinity. This last fact leads to a contradiction and the desired conclusion.

There exist quite a few related results in literature. Apart from those already
been mentioned, the work of Ni \cite{Ni} provides a classification of shrinking K\"ahler Ricci solitons with nonnegative bisectional curvature.
In a recent work of Cai \cite{Ca}, he established a classification
result for shrinking Ricci solitons under the assumption that the sectional curvature is nonnegative and bounded, and the covariant derivative of the Ricci curvature tensor decays exponentially.

\section{Proof of theorem 2}

\begin{proof}

Let us assume for the sake of contradiction that $M$ is noncompact. We begin
by recalling some important features of gradient shrinking Ricci solitons.
It is known \cite{H} that $S+\left\vert \nabla f\right\vert ^{2}-f$ is
constant on $M$, where $S$ is the scalar curvature of $M.$ So, by adding a
constant to $f$ if necessary, we may normalize the soliton such that 
\begin{equation}
S+\left\vert \nabla f\right\vert ^{2}=f.  \label{h}
\end{equation}%
Moreover, according to \cite{C}, $S>0$ unless $M$ is flat.
Concerning the potential function $f,$ Cao and Zhou \cite{CZ} have proved
that%
\begin{equation}
\left( \frac{1}{2}r\left( x\right) -c\right) ^{2}\leq f\left( x\right) \leq
\left( \frac{1}{2}r\left( x\right) +c\right) ^{2}  \label{f}
\end{equation}%
for all $r\left( x\right) \geq r_{0}$ with $c$ and $r_0$ both only depending on $n.$
Here $r\left( x\right) :=d\left(p,x\right) $ is the distance of $x$ to $p,$ a minimum point of $f$ on $M,$
which always exists. Also, the Ricci curvature $R_{ij}$ of $\left( M,g\right) $ verifies the
following differential equations. 

\begin{equation}
\Delta _{f}R_{ij}=R_{ij}-2R_{ikjl}R_{kl},  \label{eq}
\end{equation}
where $\Delta _{f}:=\Delta -\left\langle \nabla f,\right\rangle $ is the
so-called weighted Laplacian and $R_{ikjl}$ the curvature tensor. 
Let $\lambda \left( x\right) >0$ be the
smallest eigenvalue of $\mathrm{Ric}$. We note that if $v$ is an eigenvector
corresponding to $\lambda ,$ then 
\begin{equation*}
R_{ikjl}R_{kl}v_{i}v_{j}=\mathrm{Rm}\left( v,e_{k},v,e_{l}\right) R_{kl},
\end{equation*}%
where $\mathrm{Rm}$ denotes the Riemann curvature tensor. Diagonalizing the
Ricci tensor so that $R_{kl}=\lambda _{k}\delta _{kl},$ it follows from
above that 
\begin{equation}
R_{ikjl}R_{kl}v_{i}v_{j}=K\left( v,e_{l}\right) \lambda _{l}\geq 0,  \label{ps}
\end{equation}%
where $K\left( X,Y\right) $ is the sectional curvature of the plane spanned
by $X$ and $Y$. The inequality above is true because we assumed $K\left(
X,Y\right) \geq 0$, for any $X,Y.$ From (\ref{eq}) and (\ref{ps}) it follows
that $\lambda $ satisfies the following differential inequality
in the sense of barriers.
\begin{equation}
\Delta _{f}\lambda \leq \lambda .  \label{l}
\end{equation}

For a fixed geodesic ball $B_{p}\left( r_{0}\right) $ of radius $r_{0}$
large enough, let 
\begin{equation}
a:=\inf_{\partial B_{p}\left( r_{0}\right) }\lambda >0.  \label{lb}
\end{equation}

We now adapt an idea from \cite{CLY} to obtain a global lower bound for $%
\lambda $. The argument uses maximum principle applied to (\ref{l}).

Define the function 
\begin{equation}
u:=\lambda -af^{-1}-naf^{-2}.  \label{4}
\end{equation}%
By (\ref{lb}) and (\ref{4}), it follows that if $r_{0}$ is large enough
depending on dimension,%
\begin{equation}
u>0\text{ \ on \ }\partial B_{p}\left( r_{0}\right) .  \label{pos}
\end{equation}%
On $M\backslash B_{p}\left( r_{0}\right) $ we have 
\begin{eqnarray*}
\Delta _{f}\left( f^{-1}\right) &=&-\Delta _{f}\left( f\right)
f^{-2}+2\left\vert \nabla f\right\vert ^{2}f^{-3} \\
&\geq &\left( f-\frac{n}{2}\right) f^{-2}+2f^{-2} \\
&=&f^{-1}-\frac{n}{2}f^{-2},
\end{eqnarray*}%
and 
\begin{eqnarray*}
\Delta _{f}\left( f^{-2}\right) &=&2\left( f-\frac{n}{2}\right)
f^{-3}+6\left\vert \nabla f\right\vert ^{2}f^{-4} \\
&\geq &\frac{3}{2}f^{-2}.
\end{eqnarray*}%
Therefore, using this in (\ref{4}) and combining with (\ref{l}) imply
\begin{eqnarray*}
\Delta _{f}u &\leq &\lambda -af^{-1}+a\frac{n}{2}f^{-2}-\frac{3}{2}naf^{-2}
\\
&=&u.
\end{eqnarray*}%
We have thus established that 
\begin{equation}
\Delta _{f}u\leq u\text{ on }M\backslash B_{p}\left( r_{0}\right) .
\label{5}
\end{equation}

We now claim that $u\geq 0$ on $M\backslash B_{p}\left( r_{0}\right) $.
Suppose this is not true. Then there exists $x_{0}\in M\backslash
B_{p}\left( r_{0}\right) $ so that $u\left( x_{0}\right) <0$. Since by (\ref%
{pos}) we have $u>0$ on $\partial B_{p}\left( r_{0}\right) $ and $u$ is obviously nonnegative
at infinity, it
follows that $u$ achieves its minimum in the interior of $M\backslash
B_{p}\left( r_{0}\right).$ Furthermore, $u<0$ at this minimum point. By the
maximum principle, this contradicts with (\ref{5}). Thus, we conclude that $u\geq 0$ on 
$M\backslash B_{p}\left( r_{0}\right).$ So there exists some constant $0<b\leq 1$ such that 
\begin{equation}
\mathrm{Ric}\geq \frac{b}{f}\text{ \ on ~}M.  \label{6}
\end{equation}%
We use this to show that the scalar curvature  on $M\backslash
B_{p}\left( r_{0}\right) $ must satisfy

\begin{equation}
S\geq b\,\ln f.  \label{scal}
\end{equation}%
Suppose by contradiction that there exists a point $x\in M\backslash
B_{p}\left( r_{0}\right) $ where 
\[
S\left( x\right) \leq b\,\ln f\left(x\right).
\]
 Let us consider $\sigma \left( \eta \right) ,$ where $\eta \geq 0
$, to be the integral curve of $-\frac{\nabla f}{\left\vert \nabla
f\right\vert ^{2}},$ such that $\sigma \left( 0\right) =x.$ As $\left\vert
\nabla f\right\vert \left( x\right) >0$ by (\ref{h}), this flow is defined
at least in a neighborhood of $x$. Since $\frac{d}{d\eta }f\left( \sigma
\left( \eta\right) \right) =-1,$ it results that $f\left( \sigma \left( \eta
\right) \right) =t-\eta ,$ where $t:=f\left( x\right) $. Using (\ref{6}), we
have that 
\begin{eqnarray*}
\frac{d}{d\eta }S\left( \sigma \left( \eta \right) \right)  &=&-\frac{%
\left\langle \nabla S,\nabla f\right\rangle }{\left\vert \nabla f\right\vert
^{2}}=-2\mathrm{Ric}\left( \frac{\nabla f}{\left\vert \nabla f\right\vert },%
\frac{\nabla f}{\left\vert \nabla f\right\vert }\right)  \\
&\leq &-\frac{2b}{f\left( \sigma \left( \eta \right) \right) } \\
&=&-\frac{2b}{t-\eta },
\end{eqnarray*}%
where we have used the known relation that $\nabla S=2\mathrm{Ric}\left(
\nabla f\right) $. Integrating this differential inequality, we find that
\begin{equation}
S\left( x\right) -S\left( \sigma \left( \eta \right) \right) \geq 2b\,\ln
t-2b\,\ln \left( t-\eta \right) .  \label{9}
\end{equation}%
Since initially $S\left( x\right) \leq b\,\ln t$, it follows that 
\begin{eqnarray*}
S\left( \sigma \left( \eta \right) \right)  &\leq &2b\,\ln \left( t-\eta
\right)  -b\,\ln t\\
&<&b\,\ln \left( f\left( \sigma \left( \eta \right) \right) \right) .
\end{eqnarray*}
Since $S\left( y\right) \leq b\,\ln f\left( y\right) $ holds true along $%
\sigma \left( \eta \right) $ and $b\leq 1,$  (\ref{h}) implies that $\left\vert \nabla
f\right\vert \left( \sigma \left( \eta \right) \right) \geq 1$ for all  $%
0\leq \eta \leq t-1$. Hence, $\sigma \left( \eta \right) $ exists at
least for $0\leq \eta \leq t-1,$ and

\begin{equation*}
S\left( \sigma \left( t -1\right) \right) <b\,\ln f\left(\sigma (t-1)\right) =0.
\end{equation*}
This is a contradiction. In conclusion, (\ref{scal}) is true for all $x$.

Now we recall that by \cite{CZ}, for all $r>0,$
\begin{equation*}
\int_{B_{p}\left( r\right) }S\leq c\left( n\right) \mathrm{Vol}\left(
B_{p}\left( r\right) \right) .
\end{equation*}%
For any $q$ with $d(p, q)=\frac{3}{4}\,r,$ it follows from (\ref{scal}) and (\ref{f})
that 
\begin{equation}
\int_{B_{p}\left( r\right) }S\geq \int_{B_{q}\left( \frac{r}{4}\right)
}S\geq b\,\ln (\frac{r}{4}-c)^2\, \mathrm{Vol}\left( B_{q}\left( \frac{r}{4}%
\right) \right) .  \label{10}
\end{equation}%
Now Bishop-Gromov relative volume comparison implies that 
\begin{equation}
c\left( n\right) \mathrm{Vol}\left( B_{q}\left( \frac{r}{4}\right) \right)
\geq \mathrm{Vol}\left( B_{q}\left( 2r\right) \right) \geq \mathrm{Vol}%
\left( B_{p}\left( r\right) \right) .  \label{11}
\end{equation}%
By (\ref{10}) and (\ref{11}) we infer that $\ln r\leq \frac{c\left( n\right) 
}{b},$ which is a contradiction to $M$ being noncompact.

The theorem is proved.
\end{proof}

\section*{Acknowledgment}

We wish to thank Huai-Dong Cao for his interest in this work and helpful suggestions
to improve the presentation. The second author was partially supported by NSF grant No. DMS-1105799.

\end{document}